\documentclass[12pt]{amsproc}

\usepackage{amssymb}

\usepackage{graphicx}

\usepackage{amscd, amsmath, amsthm, epsf,  latexsym, color ,mathtools}

\usepackage[lmargin=1.1in,rmargin=1.1in,tmargin=1in,bmargin=1in]{geometry}
 
\usepackage{setspace} 
\usepackage{cite}
\usepackage[pdfpagelabels,bookmarksdepth=3]{hyperref}
\hypersetup{colorlinks=true,citecolor=black,urlcolor=black,linkcolor=black}
\usepackage{pinlabel}
\usepackage{pb-diagram}
\usepackage{soul}
 \usepackage{tikz-cd}
\usepackage{comment}
\usepackage{float}
\usepackage{mathrsfs}
\usepackage{wrapfig}
\restylefloat{figure}
\usepackage{caption}
\usepackage{subcaption}
\usepackage{relsize}

 \usepackage{xcolor}

\newtheorem{theoremABC}{Theorem}

\newtheorem{corollaryABC}[theoremABC]{Corollary}

\def\url#1{\expandafter\string\csname #1\endcsname}

\newtheorem{theorem}{Theorem} 
\newtheorem*{theorem*}{Theorem}

\newtheorem{proposition}[theorem]{Proposition}

\theoremstyle{definition}

\newcommand{\bbi}{{{\bf i}}}

\newcommand{\ZZ}{{\mathbb Z}}
\newcommand{\RR}{{\mathbb R}}

\newcommand{\cB}{{\B}}

\newcommand{\lto}{\longrightarrow}

\newcommand{\RN}[1]{%
  \textup{\uppercase\expandafter{\romannumeral#1}}%
}

\DeclareMathOperator{\Stab}{Stab}
\DeclareMathOperator{\Hom}{Hom}

\DeclareMathOperator{\ad}{ad}

\graphicspath{ {figures/} }

 \renewcommand{\qed}{\hfill$\square$}

\newcommand{\B}{\mathcal{B}}

\title[Homology 3-spheres whose Chern-Simons function is not Morse-Bott]{Examples of homology 3-spheres whose \\ Chern-Simons function is not Morse-Bott}

\author[H. U. Boden]{Hans U. Boden}
\address{Department of Mathematics and Statistics, McMaster University,  Hamilton, ON L8S 4K1, Canada} 
\email{boden@mcmaster.ca}

\author[C. M. Herald]{Christopher M. Herald}
\address{Department of Mathematics and Statistics, University of Nevada,  Reno, NV 89557} 
\email{herald@unr.edu}

 \author[P. Kirk]{Paul Kirk}
\address{Department of Mathematics, Indiana University, Bloomington, IN 47405} 
\email{pkirk@indiana.edu}

\thanks{HB was supported by an NSERC Discovery Grant. CH was supported by a Simons Collaboration Grant for Mathematicians. PK is thankful to Institut Fourier and the Max Planck Institute for support.}


\begin{document}

\begin{abstract}  
We construct two homology 3-spheres for which the (unperturbed) $SU(2)$ Chern-Simons function is not Morse-Bott. In one case, there is a degenerate isolated critical point. In the other, a path component of the critical set is not homeomorphic to a manifold. The examples are $+1$ surgeries on connected sums of torus knots. 
\end{abstract}
\maketitle

\section{Introduction}

The purpose of this article is to address a question raised to us by D.~ Ruberman,
 namely, whether there exist examples of {\em homology 3-spheres $M$} for which the $SU(2)$ Chern-Simons function  $$c_M \colon \cB^*\to \RR/\ZZ,$$ a circle-valued function on the space of gauge equivalence classes of irreducible $SU(2)$ connections,  fails to be Morse-Bott.  We construct an example of a homology 3-sphere whose Chern-Simons function has a degenerate isolated critical point; this example solves Problem 3.105(B) in \cite{Kirby-Problem-List}.  In addition, we construct a homology 3-sphere  for which the critical set of the  Chern-Simons function has a path component not homeomorphic to a manifold.

\medskip
 
 As is well known, holonomy identifies the critical set of $c_M$ with the irreducible character variety (a real semi-algebraic set):
 $$\chi^*(M)=\Hom(\pi_1(M),SU(2))\smallsetminus\{\theta\}/_{\rm conjugation},$$  where $\theta$ denotes the trivial homomorphism.   For any homomorphism  $\rho \colon \pi_1(M) \to SU(2)$ (henceforth called a representation), the cohomology group $H^1(M;su(2)_{\ad\rho})$ is called the {\em Zariski tangent space} of $\chi(M)$ at $\rho$.   If $M$ is a homology 3-sphere, the conjugacy class $[\theta]$ of $\theta$ is isolated in the character variety;  it follows that $\chi^*(M)$ is compact \cite{Akbulut-McCarthy}. The Hodge theorem identifies the kernel of the Hessian of $c_M$ at $\rho$ with the Zariski tangent space  of $\chi(M)$ at $\rho$ (e.g., see \cite{Taubes}). The function $c_M$ is Morse if all its critical points are non-degenerate; i.e., the Zariski tangent space is trivial at each critical point.  It is widely known that if $M$ is a connected sum of nontrivial homology spheres, $c_M$ is not Morse because $\pi_1$ is a nontrivial free product; there are gluing parameters (also known as bending parameters), related to conjugating a representation of one factor  but not the other.  
  
The Chern-Simons function $c_M:\mathcal B^* \to \RR/\ZZ$ is Morse-Bott if  every path component of the critical set  is a smooth manifold and if the Hessian of $c_M$ defines a non-degenerate bilinear form on the normal bundles of the critical submanifolds (see, for example, \cite{Nicolaescu}).  This translates into the condition that, for each $[\rho]\in \chi^*(M)$, the dimension of the Zariski tangent space of $\chi^*(M)$ at $[\rho]$ equals the dimension  of    the path component containing $[\rho]$. 

%
%
%

Fintushel-Stern \cite{Fintushel-Stern} showed that if $M$ is a Seifert-fibered homology 3-sphere, then $c_M$ is Morse-Bott. Given two homology spheres $M_1,M_2$ such that $c_{M_i}$ is Morse-Bott for $i=1,2$, the  connected sum $M_1\# M_2$  also has a Morse-Bott Chern-Simons function.  In fact, given path components $C_1\subset \chi^*(M_1)$ and $C_2\subset \chi^*(M_2)$, there are three associated components in $\chi^*(M_1\# M_2) $, diffeomorphic to $C_1\times [\theta_2]$, $[\theta_1]\times C_2$, and $C_1 \times \left( SU(2)/\{ \pm 1\} \right) \times C_2 \subset \chi^*(M_1\# M_2)$.  The latter is obtained by pairing each $\rho_1$ representing an equivalence class in $C_1$ with all $SU(2)$ conjugates of a $\rho_2$ representing a class in $C_2$. 
\medskip

Given relatively prime integers $p,q$, let $T_{p,q}$ denote the $(p,q)$ torus knot.  Consider the knot complements:  
$$X=S^3\smallsetminus nbd(T_{3,5}), Y=S^3\smallsetminus nbd(T_{2,7}), \text{ and } Z=S^3\smallsetminus nbd\left( T_{-2,7}\#T_{-2,7}\right).$$  
Equip the boundary $\partial X$ with its natural oriented meridian-longitude pair $\mu_X, \lambda_X$, and similarly $\mu_Y, \lambda_Y$ for $Y$ and  $\mu_Z, \lambda_Z$ for $Z$.  Define $h_Y:\partial X \to \partial Y$, $h_Z:\partial X \to \partial Z$ to be (orientation-reversing) homeomorphisms inducing the maps 
\begin{equation}\label{glue}
h_{Y*}\colon  \mu_X \mapsto \mu_Y, ~ \lambda_X \mapsto -\mu_Y-\lambda_Y, \hspace{.5in} 
h_{Z*}\colon  \mu_X \mapsto \mu_Z, ~ \lambda_X \mapsto -\mu_Z-\lambda_Z
\end{equation}
on the fundamental group.
Define 
$$\Sigma_1= X\cup_{h_Y} Y \ \text{ and } \ \Sigma_2=X\cup_{h_Z}  Z.$$ 
It is immediate from the fact that $X,Y,Z$ are all homology solid tori with $H_1$ generated by the meridians, and with the longitudes trivial in $H_1$, that $\Sigma_1, \Sigma_2$ are homology spheres.  

\begin{theoremABC}\label{thm1} \hfill
\begin{enumerate} 
\item There exists an isolated point in $\chi^*(\Sigma_1)$ with two-dimensional Zariski tangent space.
\item There exists a component of  $\chi^*(\Sigma_2)$  which is  not homeomorphic to a manifold.  
\end{enumerate}\end{theoremABC}

\begin{corollaryABC}
The critical set of $c_{\Sigma_1}$ contains an isolated point at which the  Hessian has  a $2$-dimensional kernel.  The critical set of $c_{\Sigma_2}$  is not homeomorphic to a manifold.\end{corollaryABC}

Thus $c_{\Sigma_1}$ and $c_{\Sigma_1}$ are most decidedly not  Morse-Bott.  Taking connected sums of these with themselves and with other homology 3-spheres provides many more complicated examples.

We note that  results of Kapovich and Millson \cite{Kapovich-Millson} imply that arbitrarily bad singularities, including isolated points with  nonzero Zariski tangent space and non-manifold path components, occur in  $SU(2)$ character varieties of 3-manifolds.  It is an open question  whether their universality results hold for  homology 3-spheres (see, e.g., \cite[Question 8.2]{Kapovich-Millson}).   We also note that there are Seifert-fibered homology spheres for which the $SU(3)$ Chern-Simons function is not Morse-Bott \cite{BHK}.  

\section{The character varieties of ${\boldsymbol X}$ and ${\boldsymbol Y}$ and their image in the character variety of the separating torus}

For any  path-connected space $A$ let 
$$\chi(A)=\Hom(\pi_1(A),SU(2))/_{\rm conjugation}$$ 
denote its character variety. Its  points are conjugacy classes, denoted $[\rho\colon \pi_1(A)\to SU(2)]$, or simply $[\rho]$.  A representation $\rho\colon \pi_1(A)\to SU(2)$  is called {\em central, (non-central) abelian}, or {\em irreducible},  depending on whether the stabilizer of $\rho$ under conjugation by $SU(2)$ is  isomorphic to $SU(2)$, $U(1)$ or $\{\pm 1\}$, respectively.

When $T$  is the 2-dimensional torus with a fixed set of generators $\mu,\lambda\in \pi_1(T)$, $\chi(T)$ is homeomorphic to a 2-sphere (usually called the pillowcase), and there is a  branched covering
\begin{equation}\label{coord}
\RR^2\to \chi(T), ~ (x,y)\mapsto [\mu\mapsto e^{x\bbi},\lambda \mapsto e^{y\bbi}]
\end{equation}
which can be seen as the composite of the projection $\RR^2\to \RR^2/(2\pi\ZZ)^2$ and the orbit map of the elliptic involution induced by $(x,y)\mapsto (-x,-y)$.  Call a curve in $\chi(T)$ a {\em line segment} if it is the image of a line segment in $\RR^2$.  Since the slope of a line is preserved by both translations by $(2\pi\ZZ)^2 $ and reflections through the origin, line segments in $\chi(T)$ have well-defined slope.

For any knot $K$, $\chi(S^3 \smallsetminus nbd(K))$ contains an  arc of (conjugacy classes of) abelian representations with central endpoints, mapping to the image of the $x$ axis (i.e., with slope zero) in $\chi(T)$.   We parameterize this arc with a path of representations  $ \mu\mapsto e^{a\bbi}, \lambda\mapsto 1 ,~ a\in [0,\pi]$, where $\mu,\lambda$ are a meridian, longitude pair.

Klassen \cite{Klassen} explicitly identified the $SU(2)$ character varieties of torus knot complements.  From his description of families of homomorphisms parameterizing the path components of $\chi^*(S^3\smallsetminus nbd(T_{p,q}))$, one can readily restrict to a meridian/longitude which generate $\pi_1(T)$  to identify the image of the restriction map 
$$ i^* \colon  \chi \left(S^3 \smallsetminus nbd(T_{p,q})\right)\to\chi(T)$$ 
induced by the inclusion $i\colon T=\partial \left( S^3 \smallsetminus nbd(T_{p,q}) \right) \to  S^3 \smallsetminus nbd(T_{p,q}).$ Along with the abelian arc, $\chi(S^3 \smallsetminus nbd(T_{p,q}))$  consists of a collection of  arcs of conjugacy classes of irreducible representations,  mapping  to $\chi(T)$ as  line segments of slope $-pq$, with ends limiting to certain points on the abelian arc.  The details in the case of $T_{3,5}$ are summarized in \cite{HHK}.  For the purposes of this article, we require only the following part of this calculation for $T_{3,5}$, $T_{2,7}$, and $T_{-2,7}$. 

\begin{proposition}[Klassen  \cite{Klassen}]  \label{Klassen result} 
 There is a path component of $\chi^*(S^3 \smallsetminus nbd(T_{3,5}))$  which is an arc mapping onto a line segment in $\chi(T)$ of slope $-15$,  $r\in \left(\tfrac{\pi}{15},\tfrac{11\pi}{15}\right)\mapsto (r, -15r)$, with  ends limiting to the points $a=\tfrac{\pi}{15}$ and $a= \tfrac{11\pi}{15}$ on the abelian arc.   Similarly, there is path component of  $\chi^*(S^3 \smallsetminus nbd(T_{\pm 2, 7}) )$   mapping onto a line segment in $\chi(T)$ of slope $\mp 14$,  $ r\in \left(\tfrac{\pi}{14},\tfrac{13\pi}{14} \right)\mapsto (r,\mp 14r)$, with ends limiting to the points $a=\tfrac{\pi}{14}$ and $a= \tfrac{13\pi}{14}$ on the abelian arc.

At each  interior point on these irreducible arcs, the Zariski tangent space is 1-dimensional.   For the (abelian) endpoints of either irreducible arc, the Zariski tangent space is 3-dimensional and  the linearization of the restriction map to $\chi(T)$ has rank one, with horizontal image. 
\end{proposition}

\begin{figure}[ht]
\begin{center}
\includegraphics[scale=1.2]{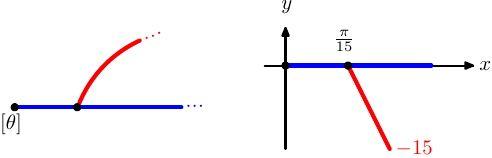}
\caption{Local picture of $\chi(X) = \chi(S^3 \smallsetminus nbd(T_{3,5}))$ near $[\theta]$ (on left) and its image under $i^*_X \colon \chi(X)\to\chi(\partial X)$ (on right) } \label{fig1}
\end{center}
\end{figure}

\begin{figure}[ht]
\begin{center}
\includegraphics[scale=1.2]{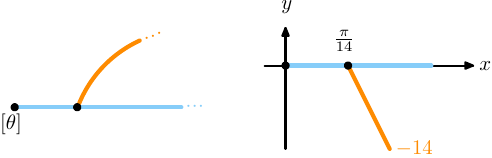}
\caption{Local picture of $\chi(Y)= \chi(S^3 \smallsetminus nbd(T_{2,7}))$ near $[\theta]$ (on left) and its image under $i^*_Y \colon \chi(Y)\to\chi(\partial Y)$ (on right) }\label{fig2}
\end{center}
\end{figure}

Figure \ref{fig1}  and Figure \ref{fig2} illustrate neighborhoods of the left ends of the irreducible arcs described in the theorem and (lifts to $\RR^2$ of) their images under restriction to the character variety of the boundary torus. In both cases, the neighborhoods embed into the pillowcase.

\begin{figure}[ht]
\begin{center}
\includegraphics[scale=1.2]{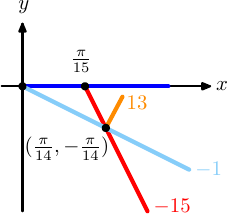}
\caption{The images  $i_X ^* \left(\chi(X)\right)$ and $h^* \circ i_Y ^* \left(\chi(Y)\right)$ near $[\theta]$ in $\chi(\partial X)$} \label{fig3}
\end{center}
\end{figure}

\section{Proof of Part (1)}  
The homeomorphism $h_Y$ of Equation (\ref{glue}) induces a map $h_Y ^*\colon \chi(\partial Y)\to \chi(\partial X)$ which lifts to the linear map 
$$h_Y ^*=\begin{pmatrix}1&0\\-1&-1\end{pmatrix}$$
on $\RR^2$, using Equation (\ref{coord}).    Figure \ref{fig3} illustrates the line segments which make up the images under the local embeddings $ i^*_X$ and $h_Y ^*\circ i^*_Y$ of the portions of $\chi(X)$ and $\chi(Y)$ in Figures \ref{fig1} and \ref{fig2}. 

Consider the fiber product
$$F:=\{ \left([\rho_X],[\rho_Y]\right) \mid i^*_X(\rho_X)=h_Y ^*\circ i^*_Y(\rho_Y)\}\subset  \chi(X)\times \chi(Y).$$
The restriction map $\chi(\Sigma_1)\to \chi(X)\times \chi(Y)$ has image $F$ and fiber over  $([\rho_X],[\rho_Y])$ (known as the space of {\em gluing parameters}) homeomorphic to the double coset space
\begin{equation}\label{double coset} 
\Stab_{\rho_X}\backslash \Stab_{\rho_{\partial X}}/\Stab_{\rho_Y}
\end{equation} 
(see, e.g., \cite{HHK}).

From the  subsets of $i_X ^* \left( \chi(X)\right), \ h_Y ^* (i_Y ^* \left( \chi(Y)\right))$ that we have identified and sketched in Figure \ref{fig3}, it is clear  that there are  two isolated points of intersection.  Specifically,  $([\theta_X],[\theta_Y])$  maps to the origin in $\chi(\partial X)=\RR^2/\sim$, and a pair $([\rho_X],[\rho_Y]) \in \chi(X)\times \chi(Y)$ which maps to $(\tfrac{\pi}{14}, -\tfrac{\pi}{14})$.  Moreover, $\rho_X$ is irreducible  and $\rho_Y$ is abelian, but non-central.

For this second pair, $\Stab_{\rho_{\partial X}}=U(1)=\Stab_{\rho_{Y}}$, since the representations $\rho_Y$ and $\rho_{\partial X}$ are abelian non-central.  It follows that $H^0(\partial X; su(2)_{ad \rho_{\partial X}})=0$. Hence the pair $([\rho_X],[\rho_Y])$ mapping to $(\tfrac{\pi}{14}, -\tfrac{\pi}{14})$ determines an {\em isolated point} of the character variety   $\chi(\Sigma_1)$, viewed as a real semi-algebraic set. 
 
Consider the Mayer-Vietoris sequence (with local $su(2)$ coefficients)
$$\cdots \stackrel{0}{\lto} H^1(\Sigma_1)\lto H^1(X)\oplus H^1(Y)\xrightarrow{i^*_X-h_Y ^* \circ i^*_Y} H^1(\partial X)\lto \cdots $$
We have:  
\begin{itemize}
\item $\dim H^1(X;su(2)_{ad\rho_X})=1$ because $[\rho_X]$ is an interior point on the irreducible arc of $\chi (X)   $ identified in Proposition \ref{Klassen result}, 
\item $\dim H^1(\partial X;su(2)_{ad\rho_{\partial X}})=2$ since the restriction $\rho_{\partial X}\colon \pi_1(\partial X)\to SU(2)$ is abelian and non-central, 
\item $\dim H^1(Y;su(2)_{ad\rho_Y})=3$ because $\rho_Y$ is the $a=\tfrac \pi{14}$ abelian endpoint of the irreducible arc of $\chi(Y)$ identified in Proposition \ref{Klassen result}, and 
\item the image of the irreducible arc in $\chi^*(X)$ and the abelian arc in $h_Y ^*(i^*_Y(\chi(Y)))$ have different slopes in $\chi(\partial X)$, namely $-15$ and $-1$, respectively,  which shows $i^*_X-h_Y ^* \circ i^*_Y$ is surjective.  
\end{itemize}
Thus, the Mayer-Vietoris sequence implies that $\dim H^1(\Sigma_1;su(2)_{ad \rho})=2$, completing the proof of the first assertion of Theorem \ref{thm1}.

\section{Proof of Part (2)} \label{Sec-4}
The proof of Part (2) of Theorem \ref{thm1} follows a similar strategy to that used to prove Part (1), but we replace $Y$ by $Z=S^3\smallsetminus nbd(T_{-2,7}\#T_{-2,7})$. The  exterior $Z$ of the composite knot $T_{-2,7}\#T_{-2,7}$ may be viewed  as the union of the two exteriors $$Z_1=Z_2=S^3\smallsetminus nbd(T_{-2,7})$$  along an annulus representing a meridian.  We begin by describing the relevant subset of $\chi(Z)$ and its image $i_Z ^* (\chi(Z)) \subset \chi(\partial Z))$.  

The fundamental group $\pi_1(Z)$ is an amalgamated free product of $\pi_1(Z_1)$ and $\pi_1(Z_2)$,  where particular meridians on each  of the two knot complements are identified.    For any representations $\rho_i:\pi_1(Z_i)\to SU(2)$, $i=1,2$, which agree on the identified meridians, there is a representation of $\pi_1(Z)$ that restricts to  $\rho_i$  on $\pi_1(Z_i)$.   The longitude for the composite knot is the product of the longitudes for $Z_1$ and $Z_2$ so, roughly speaking, the longitudinal coordinates in the pillowcase pictures for $Z_1,Z_2$ add.  
  
\begin{figure}[ht]
\begin{center}
\includegraphics[scale=1.2]{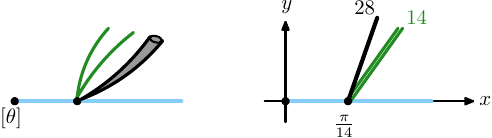}
\caption{Local picture of $\chi(Z) = \chi(S^3 \smallsetminus nbd(T_{-2,7}\# T_{-2,7}))$ near $[\theta]$ (on left) and its image under $i_Z^* \colon \chi(Z) \to \chi(\partial Z)$ (on right)} \label{fig4}
\end{center}
\end{figure}

The fiber product/gluing parameter results above (this time with restrictions to the annulus instead of to $\partial X$) demonstrate that abelian arcs and the irreducible arcs in $\chi(Z_i)$, $i=1,2$, described in Proposition \ref{Klassen result} give rise to the following subsets of $\chi(Z)$:  
\begin{enumerate}
\item[(i)] the abelian arc in $\chi(Z)$, 
\item[(ii)] two {\em half-abelian arcs}, namely an abelian/irreducible arc and an irreducible/abelian arc, consisting of representations of $\pi_1(Z)$ that are irreducible on only one of $\pi_1(Z_i)$, and 
\item[(iii)] a cylinder of irreducible/irreducible representations with $S^1$ gluing parameter.  
\end{enumerate}
All three components of $\chi(Z)$ described in (ii) and (iii) limit to the abelian points  $a=\tfrac \pi {14}, \tfrac{ 13\pi}{14}$ on the abelian arc of $\chi(Z)$.  Under $i_Z^*$, the abelian arc maps to the (image in $\chi(\partial Z)$ of) the $x$ axis; the two half abelian arcs in (ii) map to line segments of slope $14$, and the cylinder in (iii) maps onto a line segment of slope $28$.  
This is summarized in Figure \ref{fig4}.

\begin{figure}[ht]
\begin{center}
\includegraphics[scale=1.2]{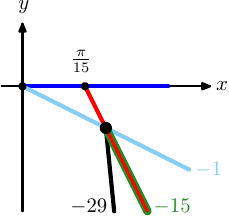}
\caption{The images  $i_X ^* \left(\chi(X)\right)$ and $h_Z ^* \circ i_Z ^* \left(\chi(Z)\right)$ near $[\theta]$ in $\chi(\partial X)$} \label{fig5}
\end{center}
\end{figure}

Under the map $h_Z ^*\colon\chi(\partial Z)\to \chi(\partial X)$, the origin is fixed (i.e., $h_Z ^*([\theta_{\partial Z}])=[\theta_{\partial X}]$), the abelian arc in $\chi(Z)$ maps to the line segment $y=-x$, the half abelian arcs described above map onto line segments leaving the abelian arc with slope $-15$, and the image of the cylinder maps onto a line segment of slope $-29$.  This is summarized in Figure \ref{fig5}.

In Figure \ref{fig5}, we overlay the images of $i_X^* (\chi(X))$ and $h_Z ^* (i_Z ^* (\chi (Z)))$ in the same picture, so that we can apply the same sort of fiber product/gluing parameter reasoning to $\Sigma_2=X\cup_{h_Z}   Z$.    We begin by noting that the abelian arc of $\chi(Z)$ meets the irreducible arc in $\chi(X)$ drawn in Figure \ref{fig1} at the point $(\tfrac \pi {14}, -\tfrac \pi {14})$.  This intersection corresponds to a point $[\rho_0]\in \chi(\Sigma_2)$ restricting to an irreducible representation  of $\pi_1(X)$ and an abelian representation of $\pi_1(Z)$, so there is no gluing parameter.  Nearby, however, the intersection includes a line segment emanating down from this point with  slope $-15$.  The preimage of that segment in $\chi(X)$ is  the irreducible arc  in Figure \ref{fig1} and the preimage in $\chi(Z)$ is the left ends of the two half-abelian arcs on in Figure \ref{fig4}.  

Taking gluing parameters into account, $[\rho_0]$ has a neighborhood in $\chi(\Sigma_2)$ which is a cone on two disjoint circles, so the path component containing $[\rho_0]$ is not a manifold.  This proves the second assertion of Theorem \ref{thm1}. \qed

\section{Further discussion and other examples}  
We note the following fact about the spaces $\Sigma_1, \Sigma_2$. 

\begin{proposition} \label{surgery on knot} The homology 3-sphere $\Sigma_1$ is diffeomorphic to $+1$ surgery on $T_{3,5}\# T_{2,7}$,  and $\Sigma_2$ is  diffeomorphic to $+1$ surgery on $T_{3,5}\# T_{-2,7}\#T_{-2,7}$.  They are both graph manifolds.  \qed
\end{proposition}

More lengthy calculations using similar techniques allow the analysis of the full character varieties of both of these homology spheres $\Sigma_i$, $i=1,2$, as well as more complicated constructions involving more torus knot complements.  We highlight a few related results without proof for the interested reader. 

\begin{proposition} 
The irreducible character variety $\chi^*(\Sigma_1)$ consists of $22$ isolated points with trivial Zariski tangent space, six isolated points with 2-dimensional Zariski tangent space like the one we described in detail, and a collection of Morse-Bott circle components.    \qed 
\end{proposition} 

While we have focused in this paper on the unperturbed Chern-Simons function, the effect of a small (carefully selected) holonomy perturbations on $\chi(\Sigma_1)$ is also reasonably straightforward to understand.  A simple holonomy perturbation in a neighborhood of $\partial X$ can be selected so that $i_X ^* (\chi(X))$ undergoes a vertical Hamiltonian flow supported away from the central endpoints on the abelian arc, so that each of the six singular isolated points resolves into a Morse critical point and the Morse-Bott circle components remain (see, for example, \cite{Herald-Kirk}).  Under a further perturbation using a curve that cuts through $\partial X$  to break the symmetry giving rise to the gluing parameters, the Chern-Simons function can be made into a Morse function; the Morse-Bott circles can be seen to each contribute two isolated critical points (contributing zero points, counted algebraically, to the Casson invariant).

For clarity, the figures only show the neighborhood of the left endpoints of the irreducible arcs described in  Proposition \ref{Klassen result}, but the irreducible arc parameterizations in that proposition show that the $T_{\pm 2, 7}$ arcs extend further to the right than the irreducible $T_{3,5}$ arc.  The following proposition is easily proved by  tracking the images of the entire half-abelian arcs (see the last two paragraphs of Section \ref{Sec-4}).
\begin{proposition} The path component of $\chi^*(\Sigma_2)$ containing the singular point $[\rho_0]$ is homeomorphic to a wedge of two 2-spheres. \qed \end{proposition}

\begin{proposition} If one replaces $Z$ with $S^3 \smallsetminus nbd(3T_{-2,7}\# T_{2,7})$ in the construction of $\Sigma_2$, then the corresponding point at $(\tfrac \pi {14}, -\tfrac \pi {14}) $ has a neighborhood that is a cone on the disjoint union of two circles and a 3-torus.  \qed\end{proposition} 

Finally, we note that the homology spheres $\Sigma_1, \Sigma_2$ can also be decomposed using Heegaard splittings, giving rise to different fiber product descriptions of the singularities described in Theorem \ref{thm1}.  For Heegaard splittings,  the character varieties of the handlebodies are smooth submanifolds of the smooth locus of the character variety of the Heegaard surface in a neighborhood of their intersections points.  Hence the local singular structure in the fiber product description associated to this decomposition is due   to this pair of handlebody character varieties intersecting nontransversely.

\bibliographystyle{alpha}

\begin{thebibliography}{HHK14}

\bibitem[AM90]{Akbulut-McCarthy}
Selman Akbulut and John~D. McCarthy.
\newblock {\em Casson's invariant for oriented homology {$3$}-spheres},
  volume~36 of {\em Mathematical Notes}.
\newblock Princeton University Press, Princeton, NJ, 1990.
\newblock An exposition.

\bibitem[BHK05]{BHK}
Hans~U. Boden, Christopher~M. Herald, and Paul~A. Kirk.
\newblock The integer valued {${\rm SU}(3)$} {C}asson invariant for {B}rieskorn
  spheres.
\newblock {\em J. Differential Geom.}, 71(1):23--83, 2005.

\bibitem[FS90]{Fintushel-Stern}
Ronald Fintushel and Ronald~J. Stern.
\newblock Instanton homology of {S}eifert fibred homology three spheres.
\newblock {\em Proc. London Math. Soc. (3)}, 61(1):109--137, 1990.

\bibitem[HHK14]{HHK}
Matthew Hedden, Christopher~M. Herald, and Paul Kirk.
\newblock The pillowcase and perturbations of traceless representations of knot
  groups.
\newblock {\em Geom. Topol.}, 18(1):211--287, 2014.

\bibitem[HK18]{Herald-Kirk}
Christopher~M. Herald and Paul Kirk.
\newblock Holonomy perturbations and regularity for traceless {$\rm SU(2)$}
  character varieties of tangles.
\newblock {\em Quantum Topol.}, 9(2):349--418, 2018.

\bibitem[Kir97]{Kirby-Problem-List}
Problems in low-dimensional topology.
\newblock In Rob Kirby, editor, {\em Geometric topology ({A}thens, {GA},
  1993)}, volume 2.2 of {\em AMS/IP Stud. Adv. Math.}, pages 35--473. Amer.
  Math. Soc., Providence, RI, 1997.

\bibitem[Kla91]{Klassen}
Eric~Paul Klassen.
\newblock Representations of knot groups in {${\rm SU}(2)$}.
\newblock {\em Trans. Amer. Math. Soc.}, 326(2):795--828, 1991.

\bibitem[KM17]{Kapovich-Millson}
Michael Kapovich and John~J. Millson.
\newblock On representation varieties of 3-manifold groups.
\newblock {\em Geom. Topol.}, 21(4):1931--1968, 2017.

\bibitem[Nic11]{Nicolaescu}
Liviu Nicolaescu.
\newblock {\em An invitation to {M}orse theory}.
\newblock Universitext. Springer, New York, second edition, 2011.

\bibitem[Tau90]{Taubes}
Clifford~Henry Taubes.
\newblock Casson's invariant and gauge theory.
\newblock {\em J. Differential Geom.}, 31(2):547--599, 1990.

\end{thebibliography}

\end{document}